\documentstyle[12pt,amssymb]{article}
\begin{document}
\newtheorem{Theo}{Theorem}
\newtheorem{Ex}{Example}
\newtheorem{Def}{Definition}
\newtheorem{Lem}{Lemma}
\newtheorem{Cor}{Corollary}
\newenvironment{pf}{{\noindent\bf Proof.\ }}{\ $\Box$\medskip}
\newtheorem{Prop}{Proposition}
\renewcommand{\theequation}{\thesection.\arabic{equation}}


\mathchardef\za="710B  
\mathchardef\zb="710C  
\mathchardef\zg="710D  
\mathchardef\zd="710E  
\mathchardef\zve="710F 
\mathchardef\zz="7110  
\mathchardef\zh="7111  
\mathchardef\zvy="7112 
\mathchardef\zi="7113  
\mathchardef\zk="7114  
\mathchardef\zl="7115  
\mathchardef\zm="7116  
\mathchardef\zn="7117  
\mathchardef\zx="7118  
\mathchardef\zp="7119  
\mathchardef\zr="711A  
\mathchardef\zs="711B  
\mathchardef\zt="711C  
\mathchardef\zu="711D  
\mathchardef\zvf="711E 
\mathchardef\zq="711F  
\mathchardef\zc="7120  
\mathchardef\zw="7121  
\mathchardef\ze="7122  
\mathchardef\zy="7123  
\mathchardef\zf="7124  
\mathchardef\zvr="7125 
\mathchardef\zvs="7126 
\mathchardef\zf="7127  
\mathchardef\zG="7000  
\mathchardef\zD="7001  
\mathchardef\zY="7002  
\mathchardef\zL="7003  
\mathchardef\zX="7004  
\mathchardef\zP="7005  
\mathchardef\zS="7006  
\mathchardef\zU="7007  
\mathchardef\zF="7008  
\mathchardef\zW="700A  

\newcommand{\be}{\begin{equation}}
\newcommand{\ee}{\end{equation}}
\newcommand{\lra}{\longrightarrow}
\newcommand{\ra}{\rightarrow}
\newcommand{\bea}{\begin{eqnarray}}
\newcommand{\eea}{\end{eqnarray}}
\newcommand{\beas}{\begin{eqnarray*}}
\newcommand{\eeas}{\end{eqnarray*}}
\newcommand{\Z}{{\Bbb Z}}
\newcommand{\R}{{\Bbb R}}
\newcommand{\C}{{\Bbb C}}
\newcommand{\SL}{SL(2,\R)}
\newcommand{\Sl}{sl(2,\C)}
\newcommand{\SU}{SU(2)}
\newcommand{\su}{su(2)}
\newcommand{\G}{{\cal G}}
\newcommand{\g}{{\frak g}}
\newcommand{\h}{{\frak h}}
\newcommand{\D}{{\rm d}}
\newcommand{\de}{\,{\stackrel{\rm def}{=}}\,}
\newcommand{\Ad}{{\rm Ad}}
\newcommand{\we}{\wedge}
\newcommand{\We}{\bigwedge}
\newcommand{\tL}{{\tilde\zL}}
\newcommand{\tl}{{\tilde\zl}}
\newcommand{\ta}{{\tilde\za}}
\newcommand{\tda}{{\widetilde{\D\za}}}
\newcommand{\nn}{\nonumber}
\newcommand{\ot}{\otimes}
\newcommand{\s}{{\textstyle *}}
\newcommand{\Li}{{\cal L}}
\newcommand{\const}{{\rm const}}
\newcommand{\pa}{\partial}
\newcommand{\ti}{\times}
\begin{center}
{\Large\bf On Filippov algebroids \\ and
multiplicative Nambu-Poisson structures}
\vskip 1cm
J. Grabowski\footnote{Institute of Mathematics, Warsaw University,
ul. Banacha 2, 02-097 Warszawa, Poland;
{\it e-mail:} jagrab@mimuw.edu.pl .\\
This work has been supported by KBN, grant No. 2 P03A 042
10.} \\ and \\
G. Marmo\footnote{Dipartimento di Scienze Fisiche, Universit\`a di Napoli,
Mostra d'Oltremare, Pad. 20, 80125 Napoli, Italy;
{\it e-mail:} gimarmo@na.infn.it .\\
This work has been partially supported by PRIN-97 "SINTESI".}
\end{center}
\date{\ }
\centerline{\bf Abstract}
We   discuss   relations   of    linear
Nambu-Poisson  structures  to   Filippov   algebras   and   define
a Filippov algebroid -- a generalization of a Lie algebroid. We also
prove   results describing multiplicative Nambu-Poisson structures
on Lie groups. In particular, it is shown that simple  Lie  groups
do not admit  multiplicative  Nambu-Poisson  structures  of  order
$>2$.
\bigskip\noindent
\setcounter{equation}{0}
\section{Introduction}
There is  growing  interest  to  $n$-ary  generalizations  of  the
concept of Lie algebra and of Poisson manifold.
In 1985 Filippov \cite{Fi} introduced a notion of $n$-Lie
algebra by assumming
that  there  is  an  $n$-linear  skew-symmetric  bracket  $$V^n\ni
(f_1,\dots,f_n)\mapsto [f_1,\dots,f_n]\in V$$ on  a  linear  space
$V$ such that the following {\em generalized Jacobi  identity}  is
satisfied:
\bea\label{FI}
&[f_1,\dots,f_{n-1},[ g_1,\dots,g_n]]=
[[f_1,\dots,f_{n-1},g_1],g_2,\dots,g_n]+\\
&[g_1,[f_1,\dots,f_{n-1},g_2],g_3,\dots,g_n]+\dots+
[g_1,\dots,g_{n-1},[ f_1,\dots,f_{n-1},g_n]].\nn
\eea
We  shall  call  such  structures  {\em   Filippov
algebras} and the identity (\ref{FI}) -- {\em Filippov identity}.
In his paper Filippov classified $n$-Lie  algebras  of  dimension
$(n+1)$  which  is  parallel  of  the  Bianchi  classification  of
3-dimensional Lie algebras.
The Filippov identity was  rediscovered  by  many  authors about
seven years later in  the  context  of   Nambu mechanics.
\par
The concept of a Nambu-Poisson structure was introduced by Takhtajan
\cite{Ta} in order to find an axiomatic formalism for
 the $n$-bracket operation
\be\label{01}
\{ f_1,\dots,f_n\}={\rm det}\left(\frac{\pa f_i}{\pa x_j}\right),
\ee
proposed  by  Nambu  \cite{Nam} (but noticed also by Filippov) to
generalize  the Hamiltonian mechanics (cf. also \cite{BF, Cha,
FDS}).
Takhtajan refers to a private communication by Flato and Fronsdal
of 1992, who observed that the Nambu canonical bracket (\ref{01})
is $n$-linear skew-symmetric and satisfies the  Filippov  identity
(he refers to it as to {\em fundamental identity}).
Such an axiom was also considered by other authors about the  same
time (see \cite{SV}).
The additional assumption made by Takhtajan was that the bracket,
acting on the algebra $C^\infty (M)$  of  smooth  functions  on  a
manifold $M$, satisfies the Leibniz rule,  i.e.
it is given by  a  multivector  field  $\zL$  ({\em  Nambu-Poisson
tensor}) on $M$ in the standard way:
\be
\{ f_1,\dots,f_n\}=\zL_{f_1,\dots,f_n},
\ee
where by $\zL_{f_1,\dots,f_k}$  we  denote  the  contraction.
$i_{\D f_k}\cdots i_{\D f_1}\zL$.
The Filippov identity means exactly that  the  Hamiltonian  vector
fields $\zL_{f_1,\dots,f_{n-1}}$ close on a Lie algebra:
\be\label{HI}
[\zL_{f_1,\dots,f_{n-1}},\zL_{g_1,\dots,g_{n-1}}]=
\sum_i\zL_{g_1,\dots,\{ f_1,\dots,f_{n-1},g_i\},\dots,g_{n-1}},
\ee
or that they preserve the tensor $\zL$, i.e. the corresponding Lie
derivatives (which we write as the Schouten bracket) vanish:
\be\label{LD0}
[\zL_{f_1,\dots,f_{n-1}},\zL]=0.
\ee
\par
The fundamental difference with the classical Poisson case is that
for   $n>2$   the  Nambu-Poisson tensor $\zL$  is decomposable
\cite{AG,Ga,MVV,Pa}, i.e. it has rank $n$ at points where it  does
not vanish. Let us note that linear tensors corresponding
to Filippov algebras  need  not  to  be  decomposable,  since  the
Filippov identity is valid not for  all  smooth  but  only for linear
functions.
\par
Similarly as we interprete elements of a Filippov algebra $V$
to be linear functions on the dual space $V^*$,
sections $Y$ of a vector bundle $\zt:E\ra M$ may be interpreted as
linear functions $\zi_Y$ on the dual bundle $\zp:E^*\ra M$:
\be
\zi_Y(\zm_m)=<Y(m),\zm_m>.
\ee
The {\em Filippov tensors} on  $E^*$  will  be,  consequently,  linear
$n$-vector fields $\zL$ on $E^*$ such that they define a  Filippov
bracket $\{\ ,\dots,\ \}_\zL$ on linear  functions  $\zi_Y$. Then,
the equation
\be\label{Fa1}
\zi_{[Y_1,\dots,Y_n]}=\{\zi_{Y_1},\dots,\zi_{Y_n}\}_\zL
\ee
defines a Filippov bracket on the space $\zG(E)$  of  sections  of
$E$. Such structures will be  called  {\em  Filippov  algebroids},
since this is a generalization of a well-known  procedure  in  the
case of linear Poisson tensors on vector bundles  and  the  notion
of Lie algebroid (cf. \cite{GU1,GU2}). Precise definitions  and
examples  will  be given in the next section.
\par
Multiplicative  Poisson  structures  are   playing   recently    a
relevant role in mathematics and physics. They can be characterized
by the  property  that  the  group  product  is  a  Poisson  map
of corresponding Poisson tensors. For multiplicative Nambu-Poisson
structures it is no longer true that the group product $G\ti  G\ra
G$ is  a  Nambu  map,  since  the  product  of  two  Nambu-Poisson
structures on $G\ti G$ is no longer a Nambu-Poisson structure.
However, for multiplicative Nambu-Poisson  structures  we  have
the corresponding infinitesimal parts being  linear  Nambu-Poisson
structures on  $\g^*$  and  the multiplicativity can be
characterized (as in the Poisson case) by the property that the
naturally defined bracket (see \ref{FF}) of (left or right)  invariant
1-forms is again an invariant 1-form (cf. \cite{Va}).
Let us also note that invariant Nambu-Poisson  structures  on  Lie
groups  are described in \cite{Na2}.
\par
There are other concepts  of  $n$-ary  Lie and Poisson
brackets using
a generalized Jacobi identity of different type than (\ref{J}) --
the skew-symmetrization of it. We will not discuss them here, so let
us only mention the papers \cite{APP,AIP,HW,ILMD,MV,SS}  and  references
there. In the recent paper \cite{VV} a unifying point of view was
proposed.
\par
The paper is organized as follows.
\par
In Section 2 we describe Filippov algebras, discuss a conjecture
stated in \cite{MVV}, and define what is  a  Filippov  algebroid
-- an $n$-ary generalization of  a  Lie  algebroid.   We   present
also  few examples of  Filippov algebroids.
\par
In Section 3 we present  some results about multiplicative
Nambu-Poisson structures, showing that they do not form as rich
family as in the classical Poisson case.
For example, we show that simple Lie  groups admit no
non-trivial multiplicative  Nambu-Poisson  structures  of
orders $>2$.  On  the  other  hand,  multiplicative  Nambu-Poisson
structures on linear spaces (regarded as commutative  Lie  groups)
are just linear Nambu-Poisson structures, i.e. a particular cases
of Filippov algebras. We get  a  description  of  all  linear
Nambu-Poisson structures similar to \cite{DZ} and \cite{MVV}.
\section{Filippov tensors and Filippov algebroids}

{\em Filippov  $n$-algebra  structures} on a vector space $V$ are
$n$-linear  skew-symmetric brackets   satisfying   the    Filippov
identity (\ref{FI}). They are determined by linear $n$-vector
fields $\zL$ on  $V^*$  which  we  shall call {\em  Filippov
tensors}.  For  a basis $x_1,\dots,x_m$ of  $V$,  regarded  as  a
basis  of  linear functions on $V^*$, we can write
$$\zL=\sum_{i_1,\dots,i_n}[x_{i_1},\dots,x_{i_n}]\pa_{x_{i_1}}\ot
\cdots\ot\pa_{x_{i_n}}.$$
In dimensions $n$  and  $(n+1)$  such  tensors  are  decomposable,
i.e. they
are  just  linear  Nambu-Poisson  tensors.  In  general,  however,
Filippov tensors may be  not  decomposable,  since,  for  example,
direct sums of Filippov  algebras  are  Filippov  algebras,  while
direct sums of non-trivial Nambu-Poisson tensors are  never
Nambu-Poisson tensors. In general, we can formulate the following.
\begin{Theo} Linear Nambu-Poisson  tensor fields of order $n>2$ are,
exactly, decomposable Filippov tensors of order $n$.
\end{Theo}
\begin{pf}
If a Filippov tensor $\zL$ of order $n$ is decomposable and
$\zL(x)\ne 0$  then (locally) it can be written as a wedge product
$\zL=X_1\we\cdots\we X_n$ for some vector fields  $X_1,\dots,X_n$.
It is aesy to see (cf. \cite{GM}, Proposition  1)  that  $\zL$  is
then Nambu-Poisson if and only if the distribution  $D$  generated
by those vector fields is involutive. Since $D$ is spanned also by
Hamiltonian vector fields $\zL_{x_{i_1},\dots,x_{i_{n-1}}}$  of  linear
functions, it is involutive in view of (\ref{HI}).
\end{pf}

In \cite{MVV} it is conjectured  that,  for
$n>2$, any $n$-Filippov algebra splits into the direct product  of
a  trivial  $n$-Filippov  algebra  and  a  number  of  non-trivial
$n$-Filippov algebras of dimensions $n$ and $(n+1)$. The following
example shows that this is not the case.

\medskip\noindent
{\bf Example 1.} Let $\zL$ be any linear Poisson tensor on the vector
space $\g^*$, corresponding therefore to a Lie  algebra  structure
on the dual $\g$. Take a basis $\{ y_1,\dots,y_m\}$  of  $\g$.  On
$V=\g^*\ti\R^k$ we   define   a   linear   $(k+2)$-vector    field
$\zL_1=\zL\we\pa_{x_1}\we\dots\we\pa_{x_k}$.  This  tensor  is   a
Filippov  tensor.  Indeed, the Filippov identity (\ref{FI})  means
that  $\zL_1$  is invariant under  the  action  of  the
Hamiltonian  vector  fields
$X=(\zL_1)_{f_1,\dots,f_{k+1}},$ where $f_i\in\{ x_1,\dots,x_k,
y_1,\dots,y_m\}$.
If we have no $y_i$ or more than two $y_i$'s among  $f_j$'s,  then
the vector field $X$ equals 0. If we have one, say $y_i$, then $X$
is proportional to $\zL_{y_i}$ and clearly
\be
[\zL_{y_i},\zL\we\pa_{x_1}\we\dots\we\pa_{x_k}]=0.
\ee
If we have two,  say  $y_i$  and  $y_j$,  then  $X$  is  a  linear
combination of $\{ y_i,y_j\}\pa_{x_s},$ where $\{ \  ,\ \}$ is the
Poisson bracket of $\zL$. Again,
\beas
&[\{ y_i,y_j\}\we\pa_{x_s},\zL\we\pa_{x_1}\we\dots\we\pa_{x_k}]=
\{ y_i,y_j\}[\pa_{x_s},\zL\we\pa_{x_1}\we\dots\we\pa_{x_k}]\pm\\
&[\zL_{y_i},\zL_{y_j}]\we\pa_{x_s}\we\pa_{x_1}\we\dots\we\pa_{x_k}
=0,
\eeas
since $\pa_{x_s}$ appears twice in the wedge product.
If $\g$ is, for example, a simple Lie algebra of  dimension  $>3$,
then the Filippov algebra structure corresponding to $\zL_1$  does
not split as in the mentioned conjecture.
\par
Indeed,  since  $\g$  is  the  derived  ideal  of   the   Filippov
algebra $V$ which does not split into ideals ($\g$ is  simple),  it
must be included in an element of the  splitting  of  $V$  into  a
direct product of ideals. But from the form of the Filippov tensor
$\zL_1$ it is clear  that  there  are  no  non-zero ideals of $V$
commuting with  the  whole  $\g$,  so  that  $V$  can  split  only
trivially and, as being a $(k+2)$-Filippov  algebra,   it   should
be  of dimension at most $(k+3)$, if the conjecture had  been
true.  The last is possible only if ${\rm dim}(\g)=3$.
Note that the Filippov tensor $\zL_1$ can be obtained by iteration
from a construction described in \cite{VV}.

\bigskip
Similarly as elements of a Filippov algebra $V$ may be interpreted
to be linear functions on the dual space $V^*$,
sections $Y$ of a vector bundle $\zt:E\ra M$ may be interpreted as
linear functions $\zi_Y$ on the dual bundle $\zp:E^*\ra M$:
\be
\zi_Y(\zm_m)=<Y(m),\zm_m>.
\ee
The {Filippov tensors} on  $E^*$  will  be,  consequently,  linear
$n$-vector fields $\zL$ on $E^*$ such that they define a  Filippov
bracket $\{\ ,\dots,\ \}_\zL$ on linear  functions  $\zi_Y$.  Now,
the equation
\be\label{Fa}
\zi_{[Y_1,\dots,Y_n]}=\{\zi_{Y_1},\dots,\zi_{Y_n}\}_\zL
\ee
defines a Filippov bracket on the space $\zG(E)$  of  sections  of
$E$. Such structures will be  called  {\em  Filippov  algebroids},
since this is a generalization of a well-known  procedure  in  the
case of linear Poisson tensors on vector bundles and Lie algebroids
(cf. \cite{GU1,GU2}).
For example, the canonical Lie algebra bracket on vector  fields,
i.e. sections of the tangent bundle $TM$, corresponds in this  way
to the linear Poisson bracket  on  the  cotangent  bundle  $T^*M$,
obtained from the canonical symplectic form. All this justifies
the following.

\medskip\noindent
{\bf Definition.} {\em A Filippov  $n$-algebroid}  is  a  vector  bundle
$\zt:E\ra M$ equipped with a Filippov $n$-bracket $[\ ,\dots,\  ]$
on sections of $E$ and a vector bundle  morphism  $a:\We^{n-1}E\ra
TM$ over identity on $M$, called {\em the anchor} of the  Filippov
algebroid, such that
\begin{description}
\item{(i)} the induced morphism on sections $a:\zG(\We^{n-1}E) \ra
\zG(TM)$ satisfies the following relation  with  respect  to   the
bracket of vector fields (cf. (\ref{HI})):
\bea
&[a(X_1\we\dots\we X_{n-1}),a(Y_1\we\dots\we Y_{n-1})]=\\
&\sum_ia(Y_1\we\dots\we [X_1,\dots,X_{n-1},Y_i]\we\dots\we Y_{n-1})\nn
\eea
\item{(ii)} and
\bea
&[X_1,\dots,X_{n-1},fY]=\\
&f[X_1,\dots,X_{n-1},Y]+a(X_1\we\dots\we X_{n-1})(f)Y\nn
\eea
for all $Y,X_1,\dots,X_{n-1}\in\zG(E)$ and $f\in C^\infty(M)$.
\end{description}

\bigskip
A study of Filippov algebroids we postpone to a separate paper. Let
us only present some examples.

\medskip\noindent
{\bf Example 2.} Let $\zL$ be a $n$-Nambu-Poisson tensor on  a  manifold
$M$ of positive dimension. The tangent (complete)  lift  $\D_T\zL$
(see  e.g.  \cite{GU1})
is never a Nambu-Poisson tensor  on  $TM$  if  $\zL\ne  0$.  In  a
coordinate system $(x^i)$ on $M$ and the adapted coordinate system
$(x^i,\dot x^j)$ on $TM$ we have
\bea\nn
&\D_T(\sum_{i_1,\dots,i_n} f_{i_1,\dots,i_n}
\pa_{x^{i_1}}\we\cdots\we\pa_{x^{i_n}})=
\sum_{k,i_1,\dots,i_n} \frac{\pa f_{i_1,\dots,i_n}}{\pa
x^k}\dot x^k\pa_{\dot x^{i_1}}\we\cdots\we\pa_{\dot x^{i_n}}+\\
&\sum_{k,i_1,\dots,i_n} f_{i_1,\dots,i_n}
\pa_{\dot   x^{i_1}}\we\cdots\we\pa_{x^{i_k}}\we\cdots\we\pa_{\dot
x^{i_n}}.
\eea
The lift $\D_T\zL$
satisfies the Filippov identity for functions
$\D_Tf=$  on  $TM$, where $f\in C^\infty(M)$ (in local coordinates
$\D_Tf=\sum_k\dot x^k\pa f/\pa x^k$). Indeed,
since $(\D_T\zL)_{\D_Tf}=\D_T(i_{\D f}\zL$) (cf.  \cite{GU1}),  we
get inductively
\be
(\D_T\zL)_{\D_Tf_1,\dots,\D_Tf_{n-1}}=\D_T(\zL_{f_1,\dots,f_{n-1}}
).
\ee
Hence,
\beas
&[(\D_T\zL)_{\D_Tf_1,\dots,\D_Tf_{n-1}},\D_T\zL]=
[\D_T(\zL_{f_1,\dots,f_{n-1}}),\D_T\zL]=\\
&\D_T[\zL_{f_1,\dots,f_{n-1}},\zL]=0,
\eeas
since the complete tangent lift  preserves  the  Schouten  bracket
(\cite{GU1}, Theorem 2.5).
\par
It is not hard to find the bracket of 1-forms induced by $\D_T\zL$:
\be\label{FF}
[\zm_1,\dots,\zm_n]=\sum_{k=1}^n(-1)^{n+k}\Li_{\zL_k}\zm_k-(n-1)\D<\zL,
\zm_1\we\cdots\we\zm_n>,
\ee
where $\Li_{\zL_k}$ denotes the Lie derivative  along  the  vector
field
\be
\zL_k=<\zL,\zm_1\we\cdots\we\check{\zm}_k\we\cdots\we\zm_n>.
\ee
This bracket can be used to determine multiplicative Nambu-Poisson
structures (cf.  \cite{Va}).
For $n>2$ this is {\bf not a Filippov algebroid} structure on $T^*M$,
since the Filippov identity is  satisfied  only for closed forms.
\medskip\noindent

{\bf Example 3.} Consider a Filippov $n$-ary bracket on $m$-dimensional
real    vector    space $V$ with    the     structure     constants
$c^k_{i_1,\dots,i_n}$ relative to a basis in  $V$.  Using  any
smooth function $g\in C^\infty(\R^m)$ we can define
a Filippov algebroid structure on the  tangent  bundle  $T\R^m$
with the trivial anchor and the $n$-ary bracket satisfying
\be
[\pa_{x^{i_1}},\dots,\pa_{x^{i_n}}]=g\sum_{k=1}^mc^k_{i_1,\dots,i_n}
\pa_{x^k}.
\ee
Explicitly,
\be
[\sum_{i=1}^mf^1_i\pa_{x^i},\dots,\sum_{i=1}^mf^n_i\pa_{x^i}]=
g\sum_{k,i_1,\dots,i_n=1}^mf^1_{i_1}\cdots f^n_{i_n}
c^k_{i_1,\dots,i_n}\pa_{x^k}
\ee
and the corresponding Filippov tensor is just
\be
\zL=g\sum_{k=1}^mc^k_{i_1,\dots,i_n}\dot x^k
\pa_{\dot x^{i_1}}\ot\cdots\ot\pa_{\dot x^{i_n}}.
\ee
\medskip\noindent
{\bf Example 4.} Consider an $(n+1)$-Filippov  algebroid  bracket  on
$T\R^m$ given by
\bea
&[\sum_{i=1}^mf^1_i\pa_{x^i},\dots,\sum_{i=1}^mf^{n+1}_i\pa_{x^i}]=\\
&\sum_{k=1}^{n+1}\sum_{i=1}^m\sum_{\zs\in
S(n)}(-1)^{k+n+1}sgn{(\zs)}f^1_{\zs(1)}\cdots f^{k-1}_{\zs(k-1)}f^{k+1}_
{\zs(k)}\cdots f^{n+1}_{\zs(n)}\frac{\pa f^k_i}{\pa x_1}\pa_{x^i},
\nn
\eea
where $S(n)$ is the group of permutations of  $(1,\dots,n)$. This
is exactly the unique Filippov algebroid structure
for which $[\pa_{x^{i_1}},\dots,\pa_{x^{i_{n+1}}}]=0$ and  the
anchor  map  is represented by the tensor field $\D
x_1\we\dots\we\D x_n\ot\pa_1$, so that the corresponding  Filippov
tensor reads
\be\zL=\pa_{\dot x^1}\we\cdots\we\pa_{\dot x^n}\we\pa_{x^1}.
\ee


\section{Multiplicative Nambu-Poisson structures}

Let  $\{\cdot,\dots,\cdot\}$  be  an   $n$-Nambu-Poisson   bracket
defined on  a  manifold  $M$.  This  means  that  the  bracket  is
$n$-linear and skew-symmetric, satisfies the Filippov identity
\bea\label{J}
&\{ f_1,\dots,f_{n-1},\{ g_1,\dots,g_n\}\}=
\{\{ f_1,\dots,f_{n-1},g_1\},g_2,\dots,g_n\}+\\
&\{ g_1,\{ f_1,\dots,f_{n-1},g_2\},g_3,\dots,g_n\}+\dots+
\{ g_1,\dots,g_{n-1},\{ f_1,\dots,f_{n-1},g_n\}\}\nn
\eea
and the Leibniz rule
\be\label{L}
\{ fg,f_2,\dots,f_{n-1}\}=f\{ g,f_2,\dots,f_{n-1}\}+
\{ f,f_2,\dots,f_{n-1}\} g.
\ee
The last means that the bracket is in fact defined by an $n$-vector
field $\zL$ in the standard way
\be
\{ f_1,\dots,f_n\}=\zL_{f_1,\dots,f_n},
\ee
where we denote $\zL_{f_1,\dots,f_k}$ to be the contraction
$i_{\D f_k}\cdots i_{\D f_1}\zL$.
The Filippov identity (\ref{J}) means then that the
{\em hamiltonian vector fields} $\zL_{f_1,\dots,f_{n-1}}$
(of $(n-1)$-tuples of functions
this time)  preserve the tensor $\zL$,
i.e. the corresponding Lie derivatives
(which we write as the Schouten bracket) vanish:
\be\label{LD}
[\zL_{f_1,\dots,f_{n-1}},\zL]=0.
\ee
This implies also that the characteristic distribution $D_\zL$ of
the $n$-vector field $\zL$, i.e. the distribution generated by
all the hamiltonian vector fields, is involutive.
Indeed, from (\ref{J}) we easily derive
\be\label{ham}
[\zL_{f_1,\dots,f_{n-1}},\zL_{g_1,\dots,g_{n-1}}]=
\sum_i\zL_{g_1,\dots,\{ f_1,\dots,f_{n-1},g_i\},\dots,g_{n-1}}.
\ee
All this look quite similar to the case of classical Poisson structures.
Now, the point is that in the case of Nambu-Poisson structures of order
$n>2$ the leaves of the characteristic foliation have to be either 0 or
$n$-dimensional,   so   that   the   Nambu-Poisson   tensor    are
decomposable.
\begin{Theo}{\rm  (\cite{AG,  Ga,  MVV,  Pa})}  If  $\zL$  is    a
Nambu-Poisson tensor of order $n>2$ not vanishing at  the  point
$p$  then  the tensor $\zL(p)$ is of rank $n$, i.e. $\zL$ is
decomposable.
\end{Theo}
We shall make later use of the following variant of the lemma `on
three planes' (cf. \cite{MVV} or \cite{DZ}).
\begin{Lem}
Let $\{\zL_i: i\in I\}$  be  a  family  of  decomposable  non-zero
$n$-vectors  of  a  vector  space  $V$   such   that   every   sum
$\zL_{i_1}+\zL_{i_2}$ is again decomposable.
Then,
\item{(a)} the linear span $D$ of the linear subspaces  $D_{\zL_i}$
they generate is at most $(n+1)$-dimensional
\item{}or
\item{(b)} the  intersection $\cap_iD_{\zL_i}$ is at least
$(n-1)$-dimensional.
\end{Lem}
\begin{pf}
It is easy to see that the sum $\zL_{i_1}+\zL_{i_2}$ is decomposable
at a point $p\in M$, where the summands are non-zero, if and only if the
intersection of $n$-dimensional subspaces
$D_{\zL_{i_1}}(p)\cap D_{\zL_{i_2}}(p)$ is at least $(n-1)$-dimensional.
Then we can use a corrected version of `lemma on three planes' as in
\cite{MVV}, Lemma 4.4., with an obvious combinatorial proof.
\end{pf}

Let now our manifold be a Lie group $G$ with the Lie algebra $\g$. It
is well known
that the tensor bundles $\We^kTG$ and $\We^kT^*G$ are canonically
Lie groups too. The group products are defined by
\be
\zL_g\circ\zL'_{g'}=(L_g)_*\zL'_{g'}+(R_{g'})_*\zL_g
\ee
and
\be
\za_g\circ\za'_{g'}=L_g^*\za'_{g'}+R_{g'}^*\za_g,
\ee
where $L_g$ and $R_g$ denote, respectively, the left and the right
translations, $(L_g)_*$ (resp. $(R_g)_*$)  are  the  corresponding
actions on contravariant tensors, and
$L_g^*$  (resp.  $R_g^*$)  is   the   dual   of   $(L_{g^{-1}})_*$
(resp. $(R_{g^{-1}})_*$).
Using    the    right    trivialization    of     the     bundles:
$\tL_g=(R_{g^{-1}})_*\zL_g$ and $\ta_g=R^*_{g^{-1}}\za_g$, we get
\be\label{m1}
\tL_g\circ\tL'_{g'}=\tL_g+\Ad_g\tL'_{g'}
\ee
and
\be
\ta_g\circ\ta'_{g'}=\ta_g+\Ad^*_g\ta'_{g'},
\ee
i.e. these groups are semidirect products  of  $G$  and  $\We^k\g$
(resp. $\We^k\g^*$), regarded as commutative groups, with  respect
to the adjoint  and coadjoint representations, respectively.
\par
A $k$-vector field $\zL:G\ra\We^kTG$ and a $k$-form $\za:G\ra\We^k
T^*G$  are  called  {\em  multiplicative}  if  they  define  group
homomorphisms.   In   other   words,   $\zL$   (resp. $\za$)   is
multiplicative if $\tL(gg')=\tL(g)+\Ad_g\tL(g')$
(resp.   $\ta(gg')=\ta(g)+\Ad^*_g\ta(g')$), i.e.    the    right
trivialization $\tL:G\ra\We^kTG$ (resp. $\ta:G\ra\We^kT^*G$) is
a 1-cocycle  of  $G$  with  coefficients  in  the  adjoint  (resp.
coadjoint) representation of $G$ in $\We^k\g$ (resp. $\We^k\g^*$).

\begin{Theo}
Let $\zL\ne 0$ be a multiplicative Nambu-Poisson structure of  order
$n>2$ on a Lie group $G$ with the Lie algebra $\g$. Then,
\item{(1)} there  is
a Lie ideal  $\h$  in  $\g$  of  dimension  $\le(n+1)$  such  that
$\tL(g)\in\We^n\h$ for all $g\in G$ (the tensor $\zL$ is therefore
tangent to the right (or left -- they are the same) cosets of  the
corresponding normal subgroup $H$ of $G$),
\item{} or
\item{(2)} there in $(n-1)$-dimensional ideal $\h$ of $\g$
such that any $\tL(g)$ is divisible by any
$0\ne\zl_0\in\We^{n-1}\h$.
\end{Theo}

\begin{pf} Denote by $V_g$ the linear subspace of $\g$ defined  by
the decomposable tensor $\tL(g)$ and by $G_0$  the  set  of  those
$g\in G$ for which $V_g\ne\{ 0\}$. Put $V_\cup={\rm span}\{ V_g:g\in
G_0\}$ and $V_\cap=\cap_{g\in G_0}V_g$.
From  $\Ad_g\tL(g')=\tL(gg')-\tL(g)$  (cf.   (\ref{m1}))   it
follows   that   $V_\cup$ and $V_\cap$  are   $\Ad_G$-invariant.
Moreover,    since for
multiplicative tensors
\be\label{x0}
\Ad_g\tL(g^{-1}g')=\tL(g')-\tL(g),
\ee
we get that  the   difference   $\tL(g')-\tL(g)$   of   two
decomposable $n$-tensors is decomposable. Hence, in view of  Lemma
1, the Lie ideal $V_\cup$ is at most  $(n+1)$-dimensional  or  the
Lie ideal $V_\cap$ is at least $(n-1)$-dimensional.
\end{pf}

Now, for any multiplicative Nambu-Poisson tesor $\Lambda$  on  $G$
let   us   define,   as   in the   Poisson case   (cf.
also   \cite{Va}),     the     corresponding      Lie      algebra
1-cocycle
$\zd_\zL:\g\ra\We^n\g$     by    the    `intrinsic     derivative'
$\zd_\zL(X)=\Li_{\tilde   X}(\tL)(e)$,
where $\tilde X$ is any vector field on $G$ with  $\tilde  X(e)=X$
(the  definition  does  not  depend  on   the   extension,   since
$\zL(e)=0$).

\begin{Theo}  The  map   $\zd^*_\zL:\We^n\g^*\ra\g^*$,   dual   to
$\zd_\zL$,  defines  a  linear   Nambu-Poisson   (in   particular,
Filippov) bracket on $\g^*$. In  other  words,  the  infinitesimal
multiplicative    Nambu-Poisson    brackets     are
`Lie-(linear-Nambu-Poisson) bialgebras'.
\end{Theo}

\begin{pf} The  linear  $n$-vector  field  $\zd_\zL$  on  $\g$  is
clearly decomposable, since, by Theorem 4, it takes values in
$\We^n\h$ and $\h$ is
$n$-  or  $(n+1)$-dimensional,  or   it   is   divisible   by   an
$(n-1)$-vector field. The generalized Jacobi  identity
for $\zd_\zL$ is a direct consequence of that for $\zL$.
\end{pf}

Now, let  us   assume   that   ${\rm   dim}(\h)=n$   in  the  case
(1) of Theorem 4 and  let  us  take
$\zl_0\in\We^n\h$,  $\zl_0\ne  0$.  Extending   $\zl_0$   by   the
right-translations to the whole $G$, we get  an  $n$-vector  field
$\zL_0$ which is tangent to the right (or left) cosets of $H$. Our
multiplicative Nambu-Poisson tensor $\zL$ is clearly  proportional
to $\zL_0$: $\zL=\zvf\zL_0$ for some smooth function $\zvf:G\ra\R$.
Let  $\zm$   be   the   modular   function   for   $\zl_0$,   i.e.
$\Ad_g\zl_0=\zm(g)\zl_0$.
The modular function is a real multiplicative  character  of  $G$,
i.e.   $\zm(gg')=\zm(g)\zm(g')$,   and   hence   (at   least   for
$G$-connected) it  is  of  the  form  $\zm=\exp(\zx)$  for  a  real
additive character $\zx:G\ra\R$: $\zx(gg')=\zx(g)+\zx(g')$.
From the multiplicativity of $\zL$ we get easily
\be
\zvf(gg')=\zvf(g)+\zm(g)\zvf(g'),
\ee
i.e. $\zvf:G\ra\R$ is a 1-cocycle with  the  coefficients  in  the
1-dimensional representation of $G$ given by $\zm$. We shall  call
such functions {\em $\zm$-characters} of $G$.
Conversely, for every  $\zm$-character  $\zvf:G\ra\R$  the  tensor
$\zvf\zL_0$ is a multiplicative Nambu-Poisson tensor on $G$.
\par
The infinitesimal part $\zd_\zvf:\g\ra\R$ satisfies
\be
\zd_\zvf([X,Y])=\zd_\zx(X)\zd_\zvf(Y)-\zd_\zx(Y)\zd_\zvf(X),
\ee
where $\zm=\exp(\zx).$ In the unimodular case  $\zm\equiv 1$,
the  Lie  algebra  1-cocycle  is   just   a   generalized   trace
\be
\zd_\zvf([X,Y])=0,
\ee
and $\zvf$ is a real (additive) character
\be
\zvf(gg')=\zvf(g)+\zvf(g').
\ee
Such characters vanish on the  derived  group  $G^{(1)}=\{  G,G\}$
and, therefore, are pull-back's of characters on  the  commutative
Lie group $G/G^{(1)}\simeq\R^k\ti T^s$ (in the connected case).
The compact part $T^s$  admits  no  additive  character,  so  that
$\zvf$ is the pull-back of a linear functional via the composition
of projections
\be
G\lra G/G^{(1)}\simeq\R^k\ti T^s\lra\R^k.
\ee
In particular, for perfect (e.g. semisimple) Lie groups we have no
non-trivial additive  characters.  On  the  other  hand,  for  the
Abelian Lie group $\R^k$ they are just linear functionals. In any
case, we have a finite-dimensional space of  additive  characters.
They can be also characterized as follows.

\begin{Theo} The following are equivalent:
\item{(a)} $f:G\ra\R$ is an additive character,
\item{(b)} $f(e)=0$ and $X^l(f)=\const$ for  every  left-invariant
vector field $X^l$ on $G$,
\item{(c)} $f(e)=0$ and $X^r(f)=\const$ for  every  right-invariant
vector field $X^r$ on $G$,
\item{(d)} $f(e)=0$ and  $\D  f$  is  a  left-and-right-invariant
1-form on $G$,
\item{(e)} $f(e)=0$ and $X^lY^l(f)=0$ (resp. $X^rY^r(f)=0$,  etc.)
for all left-invariant (resp. right-invariant, etc.) vector fields
on $G$.
\end{Theo}

\begin{pf} Let $X^l$ be the left-invariant  vector  field  on  $G$
with $X^l(e)=X\in\g$. Then,
\be
X^l(f)(g)=\frac{\D}{\D t}_{\mid t=0}f(g\exp tX)=
\frac{\D}{\D t}_{\mid t=0}(f(g)+f(\exp tX))=X(f)(e)=\const.
\ee
The rest is similar or obvious.
\end{pf}

We can also derive easily the  following result,  which  shows  that  the
polynomials  in  additive  characters can play in general the
role   of   true polynomials in the case $G=\R^k$.

\begin{Theo} The following are equivalent:
\item{(a)} $\quad f:G\ra\R$ is a polynomial of order  $\le  m$  in
additive characters on the Lie group $G$,
\item{(b)} $\quad X_1^{l(r)}\cdots  X_{m+1}^{l(r)}(f)=0\quad$  for  all
left- or right-invariant vector fields
$\quad X_1^{l(r)},\dots,X_{m+1}^{l(r)}$ on $G$.
\end{Theo}

Let us go back to the case (1) of Theorem 4. Now,
assume that ${\rm dim}(\h)=n+1$ and take  $\zl_0\in\We^{n+1}\h$,
$\zl_0\ne 0$ and the modular function  $\zm:G\ra\R$  for  $\zl_0$.
Our Nambu-Poisson tensor $\zL$  is  now  the  contraction  of  the
right-invariant prolongation $\zL_0$ of $\zl_0$ with a  tangential
1-form along the cosets of $H$, say $\za$, i.e. there  is
$\ta:G\ra\h^*$  such
that
\be\label{m2}
\tL(g)=i_{\ta(g)}\zl_0.
\ee
It is easy to see that $\ta$ is a 1-cocycle with the  coefficients
in the representation $\zm\Ad^*$ of $G$ in $\h^*$. Indeed,
\bea
\tL(gg')&=&i_{\ta(gg')}\zl_0=i_{\ta(g)}\zl_0+\Ad_gi_{\ta(g')}
\zl_0=\\
&=& i_{\ta(g)}\zl_0+i_{\Ad^*_g\ta(g')}\Ad_g\zl_0=
i_{\ta(g)}\zl_0+\zm(g)i_{\Ad_g^*\ta(g')}\zl_0,\nn
\eea
so that
\be\label{m3}
\ta(gg')=\ta(g)+\zm(g)\Ad_g^*\ta(g').
\ee
Conversely, any such $\za$ defines a  multiplicative  decomposable
tensor on $G$  by  (\ref{m2}),  since  the  values  of  $\tL$  are
automatically decomposable as $n$-tensors in an $(n+1)$-dimensional
vector  space.  The  tensor  is   a   Nambu-Poisson   tensor   if,
additionally, the   corresponding   distribution   is   involutive,
i.e.    if
$\D_H\za\we\za=0$,    where    $\D_H\za$    is    the    fiberwise
(tangential)  exterior derivative of the fiberwise (tangential)
1-form $\za$ along the  fibers--cosets of the normal subgroup $H$
of $G$.
\par
Let us note that in the case  when  $V_\za=\{\ta(g):g\in  G\}$  is
one-dimensional, we are in the previous situation, when  $\tL$  is
proportional  to  a  constant  tensor.  Indeed,  $V_\za$   is   an
$\Ad_G^*$-invariant subspace of $\h^*$, so that its annihilator in
$\h$ is  an  $n$-dimensional  Lie  ideal  in  $\g$  and  $\tL$  is
proportional to the corresponding contravariant  volume.  Thus  we
can assume in the present case that
\be\label{za}
{\rm dim}({\rm span}\{\ta(g):g\in G\})\ge 2.
\ee

Finally, let us assume that we have the case (2) of Theorem 4  and
${\rm  dim}(\h)=n-1$.  Take  $0\ne\zl_0\in\zL^{n-1}\h$.  For  each
$g\in G$ there is $X(g)\in\g$ such that
\be\label{jg}
\tL(g)=X(g)\we\zl_0
\ee
For $g\in G_0$ the vector $X(g)$ is determined  modulo  $\h$,  so,
in fact, we can regard it as a vector of $\g/\h$.  If  $\zm$  is
the modular function for $\zl_0$, then it is easy to see that
\be\label{X}
X(gg')=X(g)+\zm(g)\Ad_gX(g'),
\ee
i.e. $X:G\ra\g/\h$ is a corresponding 1-cocycle.
Conversely, any  such  cocycle  defines  a decomposable $n$-vector
field by (\ref{jg}). The corresponding distribution is  involutive
if and only if
\be\label{x}
X(g)\we\zd_X(Y)=0
\ee
for all $g\in G$ and all $Y\in\h$, where $\zd_X:\g\ra\g/\h$ is the
corresponding derived 1-cocycle. Indeed, since $\h$ is a Lie
ideal in $\g$, the corresponding distribution is spanned  by  left
(or rigt) invariant  vector  fields  on  G, corresponding  to  the
elements of $\h$, and the `vector field' $X=\sum_j\zc_jX_j$,  where
$X_j$ is a basis of $\g/\h$. The distribution is involutive  if
and only  if  the  brackets   $[Y^l,X]=\sum_jY^l(\zc_j)X_j$  ($\h$
acts trivially on $\g/\h$)  are
proportional to $X$ for $Y\in\h$. It follows from (\ref{X}) that
\be\label{quo}
[Y^l,X](g)=\zm(g)\Ad_g\Li_{Y^l}X(e)=\zm(g)\Ad_g\zd_X(Y).
\ee
The    last    one    is
proportional     to      $X(g)$      if      and      only      if
$\Ad^{-1}_gX(g)\we\zd_X(Y)=0$  for  all  $g\in  G$. As before, the
1-cocycle condition (\ref{X}) implies that $\Ad^{-1}_gX(g)$ is
proportional to  $X(g^{-1})$,  so that
\be
X(g)\we\zd_X(Y)=0
\ee
for all $g\in G,\quad Y\in\h$,
and  the  linear  span  of  values   of   $X$   is
$\Ad_G$-invariant. Hence,  $(\zd_X)_{\mid\h}\equiv 0$, or  $X$
takes values in a one-dimensional Lie ideal $\frak a$ in $\g/\h$,
spanned by the image of $(\zd_X)_{\mid\h}$. In the last case  we  are
in  the  first
situation, since our Nambu-Poisson tensor is  proportional  to  an
invariant tensor obtained from the contravariant volume on the Lie
ideal $\h'$ in $\g$, where  $\h'$  is  the inverse image of $\frak
a$ with respect to the canonical projection.
\par
In the case $(\zd_X)_{\mid\h}\equiv 0$ we can  project  $X$  to  a
mapping $\hat X:G/H\ra\g/\h$ which is the corresponding  1-cocycle
for the quotient  group  $G/H$  with  respect  to   the   quotient
of the modular function:
\be\label{hX}
\hat X([g][g'])=\hat X([g])+\hat\zm([g])\Ad_g\hat X([g']).
\ee
Indeed, from (\ref{quo}) we have $\Li_{Y^l}X(g)=0$, so that $X$ is
constant along left cosets of the corresponding Lie group $H$. The
left cosets are also right cosets ($H$ is a normal subgroup) which
implies, in view of (\ref{X}),
\be
\Li_{Y^r}X(g)=Y^r(\zm)(e)X(g),
\ee
so that $Y^r(\zm)(e)=0$ and hence $\zm$ is constant on the  cosets
of $H$ (i.e. projectable).
\par
Thus we get the following.
\begin{Theo} An $n$-vector field $\zL$ on a Lie  group  $G$  is  a
multiplicative Nambu-Poisson structure of order $n>2$ if and  only
if either
\item{(A)} there is an $n$-dimensional  normal  Lie  subgroup  $H$
of $G$  with  the  Lie  algebra  $\h$  and  the  modular  function
$\zm:G\ra\R,\quad    \zm(g)\zl_0=Ad_g\zl_0$  for  a     generator
$\zl_0\in\We^n\h$, such that
\be
\tL(g)=(R_{g^{-1}})_*\zL(g)=\zvf(g)\zl_0,
\ee
where $\zvf:G\ra\R$ is a $\zm$-character of $G$,
\item{} or
\item{(B)} there is an $(n+1)$-dimensional normal Lie subgroup $H$
of $G$  with  the  Lie  algebra  $\h$  and  the  modular  function
$\zm:G\ra\R,\quad    \zm(g)\zl_0=\Ad_g\zl_0$  for   a  generator
$\zl_0\in\We^{n+1}\h$, such that
\be
\tL(g)=(R_{g^{-1}})_*\zL(g)=i_{\ta(g)}\zl_0,
\ee
where  $\ta:G\ra\h^*$  is  a  1-cocycle satisfying (\ref{m3})  and
(\ref{za}), which defines by
$\za(g)=R^*_g\ta(g)$
a fiberwise (tangential) 1-form along the cosets of  $H$
satisfying the integrability assumption  $\D_H\za\we\za=0$,
\item{} or
\item{(C)} there is an $(n-1)$-dimensional Lie ideal $\h$ in  $\g$
and $X:G\ra\g/\h$ such that $\tL(g)=X(g)\we\zl_0$,
satisfying (\ref{X}) for the  modular  function  $\zm$  associated
with  a non-zero element $\zl_0$  of $\We^{n-1}\h$, and such  that
$X$ can be projected to $\hat X:G/H\ra\g/\h$ with (\ref{hX}).
\end{Theo}

Now, let  us  look  closer  at  the  multiplicative  Nambu-Poisson
structures in the case (B) which is the most complicated.
Denote by $\pa:\h^*\ra\h^*\we\h^*$ the Maurer-Cartan differential:
\be
<\pa\za_0,X\ot Y>=<\za_0,[X,Y]>,
\ee
$X,Y\in\h$.
\begin{Theo}
The fiberwise equation $\D_H\za\we\za=0$, in  the  presence  of
the cocycle condition (\ref{m3}),  is equivalent to the equation
\be\label{fe}
(\tda(h)-\pa\ta(g))\we(\ta(h)-\ta(g))=0
\ee
and hence to the system of equations
\bea\label{fe1}
&(a)&\qquad \D\za_{\mid H}\we\za_{\mid H}=0,\\
&(b)&\qquad
\pa\ta(g)\we\ta(g)-\pa\ta(g)\we\ta(h)-\tda(h)\we\ta(g)=0
\label{fe2}
\eea
for all $g\in G,\quad h\in H$.
\end{Theo}

\begin{pf}
Similarly to (\ref{x0}) we have  $\ta(h)-\ta(g)=\Ad^*_g\ta(g^{-1}h)$
which  shows  that  the  form   $\za_{\mid   H}-\ta(g)^r$,   where
$\ta(g)^r$   is   the   right-invariant   1-form   on   $H$   with
$\ta(g)^r(e)=\ta(g)\in\h^*$, is equivalent to the form  $\za_{\mid
g^{-1}H}$ via a diffeomorphism. Since $\D\za_{\mid g^{-1}H}\we
\za_{\mid g^{-1}H}=0$, we get (\ref{fe}) due to the  Maurer-Cartan
equation $\D(\ta(g)^r)=(\pa\ta(g))^r$.
\end{pf}

Let us consider the case $\pa\ta\we\ta\equiv 0$ in  general.  From
the  cocycle  condition   (\ref{m3})    and    the    fact    that
$V_\za={\rm span}\{ \ta(g):g\in G\}$ is $\Ad^*_G$-invariant, we
get easily that
\be\label{pa}
\pa\ta(g)\we\ta(g')+\pa\ta(g')\we\ta(g)=0
\ee
for all $g,g'\in G$.
If $\pa\ta\equiv 0$, then the $G$-coadjoint orbit of  $\ta(g)$  is
0-dimensional, i.e.  $\Ad^*_G\ta=\ta$.  It  is  not  possible  for
semisimple groups when $\ta\not\equiv 0$.
A particular case  is  $\za=\zx_1\D\zx_2$  for additive characters
$\zx_1,\zx_2$ as described above.
So let us assume that $\pa\ta\not\equiv 0$. Thus we have
$\pa\ta(g)=\ta'(g)\we\ta(g)$ for those $g\in G$ for which $\pa\ta(g)
\ne 0$ (they form an  open-dense  subset  $G_0$  of  $G$  due   to
analyticity  of multiplicative tensors) with some $\ta'(g)\in\h^*$
defined  modulo $\ta(g)$. From (\ref{pa}) it follows that
\be\label{pa1}
(\ta'(g)-\ta'(g'))\we\ta(g)\we\ta(g')=0
\ee
for all $g,g'\in G$, so that $\ta'(g')$ can be chosen from
$W_g={\rm span}\{ \ta(g),\ta'(g)\}$ for all $g\in  G_0$.
Now, we shall make use of Lemma 1 to get that ${\rm dim}V_\za\le
3$ or there  is $0\ne\za_0\in\cap_{g\in  G_0}W_g.$
It is easy to see that in the last case and with ${\rm
dim}V_\za>1$ we can write  $\pa\ta(g)=c\za_0\we\ta(g)$.
Since $V_\za$ is $\Ad^*_G$-invariant and  $\Ad^*_g\pa=\pa\Ad^*_g$,
the covector $\za_0$ is $\Ad^*_G$-invariant. Again, this is
impossible in the semisimple case.
\par
Now, let us consider the case  $\pa\ta(g)\we\ta(g)\ne  0$  for  $g$
from an open-dense subset $G_0$ of $G$. Together with  (\ref{fe2})
it gives that $\tda(e)\ne 0$ and
\be
\tda(e)=\za_1\we\za_2\ne 0
\ee
for some $\za_1,\za_2\in\h^*$ , since by  (\ref{fe1})  $\tda(h)$
is decomposable for $h$ from an open-dense subset of $H$,  thus
for all $h\in H$.
Using multiplicativity to
\be
\pa\ta(g)\we\ta(g)-\za_1\we\za_2\we\ta(g)=0,
\ee
which is (\ref{fe2}) for $h=e$, we get
\bea\label{fe4}
&\pa\ta(g)\we\Ad^*_g\ta(g')+\Ad^*_g\pa\ta(g')\we\ta(g)+\\
&\zm(g)\Ad^*_g\pa\ta(g')\we\Ad^*_g\ta(g')-\za_1\we\za_2\we\Ad^*_g
\ta(g')=0.\nn
\eea
Taking the wedge product of (\ref{fe4}) with $\Ad^*_g\ta(g')$ and
using the $\Ad^*_G$-invariance of $V_\za$, we get finally
\be
\pa\ta(g')\we\ta(g')\we\ta(g)=0
\ee
for all $g,g'\in G$. Hence  ${\rm   rank}(\pa\ta(g))=2$   and  ${\rm
dim}V_\za\le 3$. We  can summarize our considerations as follows.

\begin{Theo} If the multiplicative Nambu-Poisson structure  is  of
type  (B),  Theorem  5,   then   the   1-cocycle    $\ta:G\ra\h^*$
satisfies equations (a) and (b) of Theorem 8 and takes
values in $G$-coadjoint orbits of dimension $\le  2$,  i.e.  ${\rm
rank}(\pa\ta(g))\le 2$. If ${\rm rank}(\pa\ta(g))=2$ for  some  $g\in
G$,   then   $V_\za={\rm    span}\{\ta(g):g\in    G\}$    is    an
$\Ad^*_G$-invariant subspace of $\h^*$ of dimension $\le 3$, or
there  is  an   $\Ad^*_G$-invariant   $\za_0\in\h^*$   such   that
$\pa\ta(g)=\za_0\we\ta(g)$.
\end{Theo}
\begin{Cor} For a semisimple Lie group $G$ with the  decomposition
into   simple   factors   $G=G_1\ti\cdots\ti   G_m$    the    only
multiplicative Nambu-Poisson tensors $\zL$  of  order   $>2$   are
wedge products of the contravariant volume on a part, say  $G_1\ti
\cdots\ti G_k$, of the decomposition with
\item{(1)} either a multiplicative Lie-Poisson tensor on a
3-dimensional factor (so that $\SL$ or $\SU$ must appear in the
decomposition),
\item{(2)}  or  a  multiplicative  vector  field   on   the   rest
$G_{k+1}\ti\cdots\ti G_m$ of the decomposition (which is always the
difference of the left and right prolongation of an  element  from
the corresponding Lie algebra).
\par
In particular, simple  Lie  groups  do  not  admit  multiplicative
Nambu-Poisson structures of order $> 2$.
\end{Cor}
\begin{pf} Since semisimple  Lie  groups  do  not  admit  additive
characters, the case (A) of Theorem 5 is not  possible. In  the
case (C) it is easy to see that $X:G\ra\g/\h$ is projectable to  a
multiplicative vector field on $G/H\simeq G_{k+1}\ti\cdots\ti G_m$.
Let us consider the case (B).
In the semisimple case the normal subgroup $H$ is a  part  of  the
decomposition, say $H=G_1\ti\cdots\ti G_k$.  Since  the  coadjoint
action of $G$ on $\h^*$ reduces to the  coadjoint  action  of  $H$
(the rest of simple factors acts  trivially)  and  the  orbits  of
non-zero elements in  $\h^*$  are  of   dimension   $\ge   2$,  we
conclude that ${\rm dim}V_\za\le 3$ and the annihilator of $V_\za$
is a Lie  ideal  in
$\g$  of  codimension  $\ge  3$,  thus  $=3$.  Hence  there  is  a
3-dimensional simple factor, say $G_1$ and $\ta$ takes  values  in
$\g_1^*$. Moreover,
\be
\tL(g)=i_{\ta(g_1)}\tl_1(g_1)\we\tl_1(g_2)\we\cdots\we\tl_k(g_k),
\ee
where $0\ne\tl_i\in\We^{{\rm dim}\g_i}\g_i$.
\end{pf}

The other extreme case are commutative  groups  $G=\R^k$.  We  get
easily the following generalization of the $n$-Bianchi
classification  of Filippov brackets in \cite{MVV} (cf. also
\cite{DZ}).
\begin{Cor} Every linear Nambu-Poisson  tensor of order $n>2$ on
$\R^m$  is  in suitable linear coordinates of the form
\item{(A)}
\be
\zL=\zvf\pa_1\we\cdots\we\pa_n,
\ee
where $\zvf$ is any linear function, or
\item{(C)}
\be
\zL=\pa_1\we\cdots\we\pa_{n-1}\we(\sum_{i,j\ge n}a_{ij}x_j\pa_j),
\ee
for some constants $a_{ij}$, or
\item{(B)}
\be
\zL=i_\za(\pa_1\we\cdots\we\pa_{n+1}),
\ee
where the linear 1-form $\za$ is characterized as follows:
\item{(1)}
in the case
$\D(\za_{\mid\R^{n+1}})=0$,
\be
\za=\D \zvf+\sum_{i>n+1\atop j\le n+1}a_{ij}x_i\D x_j,
\ee
where   $\zvf$   is   a   quadratic   polynomial   in    variables
$x_1,\dots,x_{n+1}$ and $a_{ij}\in\R$;
\item{(2)}
in the case
$\D(\za_{\mid\R^{n+1}})(0)\ne 0$,
\be
\za=\D \zvf+\frac{1}{2}((x_1+\sum_{i>2}a_ix_i)\D x_2
-(x_2+\sum_{j>2}b_jx_j)\D x_1),
\ee
where $\zvf$ is a quadratic polynomial in variables $x_1,x_2$  and
$a_i,b_j\in\R$.
\end{Cor}

\begin{pf}  The  cases  (A)  and  (C)  follow  easily   from   the
corresponding cases of Theorem 4. For
the case (B) and the commutative group, $\pa\ta\equiv  0$ and  we  may
identify $\za$ with $\ta$. The normal subgroup $H$ is  clearly  an
$(n+1)$-dimensional  subspace,   the   contravariant   volume   is
unimodular, and the cocycle condition for $\ta$ is just  linearity
of  the  1-form  $\za=\sum_{j\le   n+1}a_{ij}x_j\D   x_j$.   If
$\D(\za_{\mid\R^{n+1}})(0)=0$, we have no  more   restrictions  in
view of (\ref{fe}).  In
the  other   case,   $\D(\za_{\mid\R^{n+1}})=\const\ne   0$.   This
2-covector is decomposable, say, equal to $\D x_1\we\D  x_2$.
Then $\za$ takes only values spanned by $\D x_1$ and $\D  x_2$  in
view of (\ref{fe2}).
\end{pf}

\medskip\noindent
{\bf Example 5.}
Let $\zP$ be the Poincar\'e group. Topologically, it  is  a  direct
product $\zP=L\ti\R^4$ of the Lorenz group $L$ and $\R^4$ with the
semidirect  group  multiplication  $(g,x)\circ(g',x')=(gg',x+gx')$
relative to the canonical  action  of  $L$  on  the  Minkowski  space
$\R^4$.    This    action    preserves    the    Minkowski    form
$\zvf=x_0^2-\sum_{i=1}^3x_i^2$, so that  its  exterior  derivative
$\za=(1/2)\D\zvf$ is a multiplicative 1-form on $\zP$ which is
tangent  to  the cosets of $\R^4$. The tensor $\zL_0=\pa_{x_0}
\we\pa_{x_1}\we\pa_{x_2}\we\pa_{x_3}$,
representing   the   contravariant    volume    on    $\R^4$    is
$\zP$-invariant, so that the contraction $\zL=i_{\za}\zL_0$  is  a
multiplicative 3-vector field on $\zP$. It is clearly  involutive,
since $\za$ is a closed form. In coordinates,
\beas
\zL&=x_0\pa_{x_1}\we\pa_{x_2}\we\pa_{x_3}+
x_1\pa_{x_0}\we\pa_{x_2}\we\pa_{x_3}-\\
&x_2\pa_{x_0}\we\pa_{x_1}\we\pa_{x_3}+
x_3\pa_{x_0}\we\pa_{x_1}\we\pa_{x_2}.
\eeas
The  hamiltonian  vector   fields   of   pairs   $(x_0,x_i)$   and
$(x_i,x_j)$, $i,j=1,2,3$ are
\beas
X_{x_0,x_1}=x_3\pa_{x_2}-x_2\pa_{x_3},
& X_{x_0,x_2}=x_1\pa_{x_3}-x_3\pa_{x_1}, &
X_{x_0,x_3}=x_2\pa_{x_1}-x_1\pa_{x_2},\\
X_{x_1,x_2}=x_0\pa_{x_3}+x_3\pa_{x_0}, &
X_{x_3,x_1}=x_0\pa_{x_2}+x_2\pa_{x_0}, &
X_{x_2,x_3}=x_0\pa_{x_3}+x_3\pa_{x_0}.
\eeas
You can recognize the usual duality between  the  plane  in  which
takes place the (pseudo) rotation and the plane  `normal'  to  it,
like in three dimensions  between  the  plane  and  the  `axis  of
rotation' normal to it.

\bigskip\noindent
{\bf Acknowledgment} We are grateful to Jean Paul Dufour  for  his
remarks on the formulation of Lemma 1 and to Juan  Carlos  Marrero
for remarks on Filippov algebroids.

\end{document}